\documentclass{amsart}

\usepackage{amsfonts,amssymb,amscd,amstext}
\usepackage[a4paper,margin=4cm]{geometry}
\usepackage{hyperref}
\usepackage{esint}

\usepackage{graphicx}

\renewcommand{\leq}{\leqslant}

\renewcommand{\le}{\leqslant}
\renewcommand{\ge}{\geqslant}
\newcommand{\ptl}{\partial}

\newcommand{\rr}{{\mathbb{R}}}
\newcommand{\cc}{{\mathbb{C}}}
\newcommand{\bb}{{\mathbb{B}}}
\newcommand{\la}{\lambda}

\newcommand{\hh}{{\mathbb{H}}}

\newcommand{\sph}{{\mathbb{S}}}

\newcommand{\hhn}{\mathbb{H}^n}
\newcommand{\sub}{\subset}

\newcommand{\escpr}[1]{\big<#1\big>}
\newcommand{\Sg}{\Sigma}
\newcommand{\Om}{\Omega}
\newcommand{\eps}{\varepsilon}

\newcommand{\ga}{\gamma}

\newcommand{\nuh}{\nu_{H}}

\newcommand{\area}{A}
\newcommand{\vol}{V}

\DeclareMathOperator{\divv}{div}

\DeclareMathOperator{\intt}{int}

\newtheorem{theorem}{Theorem}[section]

\newtheorem{lemma}[theorem]{Lemma}

\theoremstyle{definition}

\theoremstyle{remark}

\numberwithin{equation}{section}

\begin{document}

\title[An isoperimetric inequality in the Heisenberg group $\hh^n$]{A
proof by calibration of an isoperimetric inequality in the Heisenberg
group $\hh^n$}

\author[M.~Ritor\'e]{Manuel Ritor\'e} \address{Departamento de
Geometr\'{\i}a y Topolog\'{\i}a \\
Universidad de Granada \\ E--18071 Granada \\ Espa\~na}
\email{ritore@ugr.es}

\date{\today}

\thanks{Research supported by MCyT-Feder grant MTM2010-21206-C02-01 and Junta de Andaluc\'{\i}a grants FQM-325 and P09-FQM-5088.} 
\subjclass[2000]{53C17, 53C42, 49Q20}
\keywords{Sub-Riemannian geometry, Heisenberg group, area-stationary surface, constant mean curvature surface, isoperimetric problem, isoperimetric region, calibration}

\begin{abstract}
Let $D$ be a closed disk centered at the origin in the horizontal hyperplane $\{t=0\}$ of the sub-Riemannian Heisenberg group $\hh^n$, and $C$ the vertical cylinder over $D$.  We prove that any finite perimeter set $E$ such that $D\subset E\subset C$ has perimeter larger than or equal to the one of the rotationally symmetric sphere with constant mean curvature of the same volume, and that equality holds only for these spheres using a recent result by Monti and Vittone \cite{mv}.
%If the boundary of $E$ is locally lipschitz in Euclidean sense, or an $\hh$-regular hypersurface, then equality holds for the rotationally symmetric spheres.
\end{abstract}

\maketitle

\thispagestyle{empty}

\bibliographystyle{amsplain} %\nocite{*}

\section{Introduction}
\label{sec:intro}
\setcounter{equation}{0}

It was conjectured by P.~Pansu in 1983 \cite{pansu1} that the isoperimetric regions, minimizing perimeter under a volume constraint, in the sub-Riemannian Heisenberg group $\hh^1$ are the topological balls enclosed by the one-parameter family $\{\sph_\la\}_{\la>0}$ of rotationally symmetric spheres of constant mean curvature $\la$ described in \cite{pansu1}, see also \cite{leomas}, \cite{rr1}.

The existence of isoperimetric regions in Carnot groups was proven by G.-P.~Leonar\-di and S.~Rigot \cite{lr}. Every Carnot group $\mathbb{G}$ is equipped with an one-parameter family of dilations which has a well-known effect on the perimeter and the Haar measure, the volume of the group. Hence the isoperimetric profile of $\mathbb{G}$, the function assigning to each volume $v>0$ the infimum of the perimeter of the sets of volume $v$, is given by
\[
I_\mathbb{G}(v)=C_\mathbb{G}\,v^{(Q-1)/Q},
\]
where $Q>0$ is the homogeneous dimension of $\mathbb{G}$, and $C_\mathbb{G}>0$ is a constant. The isoperimetric profile must be seen as an optimal isoperimetric inequality in $\mathbb{G}$. For any finite perimeter set $E\subset\mathbb{G}$,
\begin{equation}
\label{eq:isopg}
\tag{*}
P_\mathbb{G}(E)\ge C_\mathbb{G}\,|E|^{(Q-1)/Q},
\end{equation}
where $P_\mathbb{G}$ is the sub-Riemannian perimeter and $|E|$ is the volume of $E$. Pansu's conjecture then states that, for $\mathbb{G}=\hh^1$, which has homogeneous dimension $Q=4$, equality is attained in \eqref{eq:isopg} precisely when $E$ is the topological ball enclosed by some sphere $\sph_\la$. This conjecture can be extended to the higher dimensional Heisenberg groups $\hh^n$, $n\ge 2$.

Several attempts to solve this conjecture have been made. R.~Monti \cite{monti} and G.-P.~Leonardi and S.~Masnou \cite{leomas} have shown that there is no direct counterpart in $\hh^1$ to the Brunn-Minkowski inequality in Euclidean space.  In fact, such a Brunn-Minkowski type inequality would imply that the metric balls for the Carnot-Carath\'edory distance would be isoperimetric regions, which is known to be false \cite{monti2}. The author and C.~Rosales showed in \cite{rr1} that the only compact rotationally symmetric $C^2$ hypersurfaces in $\hh^n$ with constant mean curvature are the spheres $\sph_\la$, see also \cite{ni}. The same authors proved in \cite{rr2} that the spheres $\sph_{\la}$ are the only compact $C^2$ surfaces in $\hh^1$ which are area-stationary under a volume constraint, thus solving the isoperimetric problem assuming $C^2$ regularity of the solutions.  R.~Monti and M.~Rickly \cite{mm} proved that the spheres $\sph_{\la}$ are isoperimetric in $\mathbb{H}^1$ under the additional assumption of Euclidean convexity.   In \cite{dgn06}, D.~Danielli et al. provided a proof of the isoperimetric property of the spheres $\sph_{\la}\subset\hh^n$ in the class of sets that are the union of the graphs of a non-negative and a non-positive function and a negative graph of class $C^2$ over a Euclidean disk centered at the origin in the horizontal hyperplane $\{t=0\}$ in $\hh^n$, and
enclosing the same volume above and below such hyperplane. For a description of these results, and some other approaches, the reader may consult Chapter~8 of the monograph by L.~Capogna et al. \cite{cdpt}.

In this paper we extend the main result in \cite{dgn06}. We prove in Theorem~\ref{th:main} that, if $C$ is the vertical cylinder in $\hh^n$ over a closed disk $D$ centered at the origin in the horizontal hyperplane $\{t=0\}$, and $E\subset\hh^n$ is a finite perimeter set so that $D\subset E\subset C$, then the perimeter of $E$ is larger than or equal to the one of the ball $\mathbb{B}_\la$ enclosed by the sphere $\sph_\la$ with $|\bb_\la|=|E|$. Equality characterizes the spheres $\sph_\la$ by a recent result of R.~Monti and D.~Vittone \cite{mv}, who proved that a set in $\hh^n$ of locally finite perimeter with continuous horizontal unit normal has $\hh$-regular boundary. Theorem~\ref{th:main} can be applied to a set $E\subset\hh^n$ rotationally symmetric with respect to a vertical axis passing through the origin and satisfying $D\subset E$. Assumption $D\subset E$ has been recently removed by R.~Monti \cite{monti3}, who has proven, using Theorem~\ref{th:main} and a symmetrization argument, that the spheres $\sph_{\la}\subset\hh^n$ are isoperimetric in the class of rotationally~symmetric sets of finite perimeter. Theorem~\ref{th:main} and the result of R.~Monti and D.~Vittone \cite{mv} together provide the only known characterization result for solutions of the isoperimetric problem in $\hh^n$ in the class of finite perimeter sets.

One of the classical proofs of the isoperimetric inequality in $\rr^3$ was given by H.~A.~Schwarz \cite{schwarz}, and extended to higher dimensional Euclidean spaces, spheres and hyperbolic spaces in a series of papers by E.~Schmidt \cite{schrn}, \cite{sch}, \cite{MR0002196}, \cite{MR0009127}. Schwarz's proof had two main ingredients: what it is known nowadays as Schwarz's symmetrization and a variational argument to prove that in the class of rotationally symmetric sets spheres have the smallest perimeter under a volume constraint. A unified argument, using calibrations, can be used to give a proof of the second part of Schwarz's argument in a wide class of homogeneous Riemannian manifolds \cite[\S~1.3.1]{ritsin}, and in sub-Riemannian manifolds. A symmetrization in the sub-Riemannian Heisenberg group $\hh^n$ is difficult to produce due to the lack of reflections with respect to hyperplanes, on which all classical symmetrizations are based.

In the proof of Theorem~\ref{th:main} we consider the right cylinder $C\subset\hh^n$ over a closed Euclidean disk $D$ in the horizontal hyperplane $\{t=0\}$. On $C$ we construct two foliations by vertically translating the upper hemisphere and the lower one of the only sphere $\sph_\la$ intersecting $\{t=0\}$ at $\ptl D$. Using these foliations we prove that the sphere $\sph_\la$ minimize the functional $\text{area}-n\la\,\text{volume}$ in the class of finite perimeter sets $E\subset\hh^n$ satisfying $D\subset E\subset C$. Then we minimize over the spheres $\sph_\mu$ the functional $\text{area}-n\mu\,(\text{volume}-|E|)$ to get the desired result. The reader should compare our proof with the one given by E.~Schmidt \cite{sch} of the isoperimetric property of balls in the $n$-dimensional sphere $\sph^n$ in the class of rotationally symmetric sets. 

%in Theorem~\ref{th:main} to any set $E$ with finite perimeter contained in the
%cylinder $C$ over a disk $D$ centered at the origin in the $t=0$
%hyperplane and containing the disk.  In our proof a calibration
%argument is used to compare the perimeter of $E$ and its volume.  with
%the one of the sphere $\sph_{\la}$ inscribed in the cylinder.  Then
%the perimeter of $E$ is compared with the one of the sphere enclosing
%volume $|E|$.  The equality case is analyzed by considering the
%characteristic curves in $\ptl E$.  The reader should compare our
%proof with the one by E.~Schmidt \cite{sch} of the isoperimetric
%property of balls in the $n$-dimensional sphere $\sph^n$ in the class
%of rotationally symmetric sets.  The proof is valid in a wide class of
%Riemannian manifolds, including simply connected space forms and
%homogeneous manifolds.  The results in this paper have been recently
%used by R.~Monti \cite{monti3} to prove that the spheres $\sph_{\la}$
%are isoperimetric in the class of rotationally symmetric sets.

We have organized this paper into two sections.  In the following one we state some material needed in the proof of Theorem~\ref{th:main}. Proofs of the results which are essential but cannot  be found in the literature, in particular of Lemmae~\ref{lem:a-nhv}, \ref{lem:char}, \ref{lem:spheres} are outlined. In section~\ref{sec:isop} we give the proof of our main result Theorem~\ref{th:main}.

The author is extremely grateful to Roberto Monti and Davide Vittone for sending him a copy of their manuscript \cite{mv}, and to the referee for his valuable suggestions.

\section{Preliminaries}
\label{sec:preliminaries}

\subsection{The Heisenberg group}

The Heisenberg group $\hh^n$ is the Lie group $(\rr^{2n+1},\cdot)$, where we consider in $\rr^{2n+1}\equiv\cc^n\times\rr$ its usual differentiable structure and the product
\[
(z,t)\cdot(w,s)=(z+w, t+s+\sum_{i=1}^n \text{Im}(z_i\bar{w}_i)).
\]
A basis of left-invariant vector fields is given by $\{X_1,...,X_n,Y_1,...,Y_n,T\}$, where
\[
X_i=\frac{\ptl}{\ptl x_i}+y_i\,\frac{\ptl}{\ptl t},\quad
Y_i=\frac{\ptl}{\ptl y_i}-x_i\,\frac{\ptl}{\ptl t},\quad i=1,\dots, n;
\qquad
T=\frac{\ptl}{\ptl t}.
\]
The only non-trivial bracket relations are $[X_i,Y_i]=-2T$, $i=1,\ldots,n$.  The \emph{horizontal distribution} at a point $p\in\hh^n$ is defined by $\mathcal{H}_{p}:=\text{span}\{(X_{i})_{p}, (Y_{i})_{p}: i=1,\ldots,n\}$.

%\subsection{The left-invariant metric}

We shall consider on $\hh^n$ the left-invariant Riemannian metric $g=\escpr{\cdot,\cdot}$ so that the basis $\{X_1,...,X_n,Y_1,...,Y_n,T\}$ is orthonormal.  The \emph{horizontal
projection} of a vector field $U$ in $\hh^n$, denoted by $U_{H}$, is the orthogonal projection of $U$ over $\mathcal{H}$.

The Levi-Civita connection on $(\hhn,g)$ is denoted by $D$. From Koszul formula and the Lie bracket relations we get
\begin{alignat}{2}
\notag D_{X_i}X_j&=D_{Y_i}Y_j=D_{T}T=0, \\
\label{eq:christoffel}
D_{X_i}Y_j&=-\delta_{ij}\,T, \qquad D_{X_i}T=Y_i, \qquad D_{Y_i}T=-X_i, \\
\notag D_{Y_i}X_j&=\delta_{ij}\,T, \qquad \ \ \,D_{T}X_i=Y_i, \qquad
D_{T}Y_i=-X_i.
\end{alignat}
For any vector field $U$ on $\hhn$ we define $J(U):=D_UT$.  It follows
from \eqref{eq:christoffel} that $J(X_i)=Y_i$, $J(Y_i)=-X_i$, and
$J(T)=0$, so that $J$ defines a linear isometry when restricted
to the horizontal distribution.
% Note also that $\escpr{J(U),V}+\escpr{U,J(V)}=0$, for any pair of
% vector fields $U$ and $V$.

\subsection{Volume and sub-Riemannian perimeter} The \emph{volume} $|E|$ of a Borel set $E\subset\hh^n$ is the Riemannian volume of $E$ with respect to the metric $g$.  The \emph{area} $A(\Sg)$ of a $C^1$ hypersurface $\Sg\subset\hh^n$ is defined as
\[
A(\Sg):=\int_{\Sg}|N_{H}|\,d\Sg,
\]
where $d\Sg$ is the area element induced on $\Sg$ by the Riemannian metric $g$, and $N$ is a locally defined unit vector normal to $\Sg$.  For a $C^2$ hypersurface enclosing a bounded region $E$, the area coincides with the sub-Riemannian perimeter $|\ptl E|$,  defined as
\[
|\ptl E|(\Om):=\sup\bigg\{\int_{\Om}\divv U\,dv : U \text{ 
horizontal of class}\ C^1, |U|\le 1, \text{supp}(U)\subset\Om
\bigg\},
\]
where $\Om\subset\hh^n$ is an open set, $\divv U$ is the Riemannian divergence of the vector field $U$, and $dv$ is the volume element associated to $g$.

A set $E\subset\hh^n$ is of \emph{locally finite perimeter} if $|\ptl E|(\Om)<+\infty$ for all bounded open sets $\Om\subset\hh^n$.  It is of \emph{finite perimeter} if $|\ptl E|:=|\ptl E|(\hh^n)<+\infty$.  We normalize any finite perimeter set $E\subset\hh^n$ to include its density one points and to exclude its density zero points \cite[Prop.~3.1]{MR775682}.

The \emph{reduced boundary} $\ptl^*E$ of a set $E\subset\hh^n$ of locally finite perimeter, as defined in \cite[Def.~2.17]{fssc}, is composed of those points $p\in\hh^n$ such that
\begin{enumerate}
\item $|\ptl E|(U(p,r))>0$ for all $r>0$,
\item there exists $\lim_{r\to 0} \fint_{U(p,r)}\nu_H\,d|\ptl E|$,
\item $\big|\lim_{r\to 0} \fint_{U(p,r)}\nu_H\,d|\ptl E|\big|=1$.
\end{enumerate}
Here $U(p,r)$ is the open metric ball with respect to the distance induced by the homogeneous norm
\begin{equation}
\label{eq:norm}
||p||_\infty:=\max\{|z|,|t|^{1/2}\}, \qquad p=(z,t),
\end{equation}
which is globally equivalent to the Carnot-Carath\'eodory distance, \cite[Prop.2.7]{fssc}. The perimeter measure of a set of locally finite perimeter is supported on the reduced boundary as shown in \cite[Thm.~7.1]{fssc}.

\subsection{Hypersurfaces in $\hh^n$ and variational formulae}

For a $C^1$ hypersurface $\Sg\sub\hhn$, the \emph{singular set} $\Sg_0\subset\Sg$ consists of the points where the tangent hyperplane coincides with the horizontal distribution.  The set $\Sg_0$ is closed and has empty interior in $\Sg$, and so the \emph{regular set} $\Sg\setminus\Sg_0$ of $\Sg$ is open and dense in $\Sg$.  For any $p\in\Sg\setminus\Sg_0$, the tangent hyperplane meets transversally the horizontal distribution, and so $T_{p}\Sg\cap\mathcal{H}_{p}$ is $(2n-1)$-dimensional. We say that $\Sg$ is \emph{two-sided} if there is a globally defined unit vector field normal to $\Sg$. Every $C^1$ hypersurface is locally two-sided.

Let $\Sg$ be a $C^2$ hypersurface in $\hhn$, and $N$ a unit vector
normal to $\Sg$.  The singular set $\Sg_0\sub\Sg$ can be described as
$\Sg_{0}=\{p\in\Sg : N_H(p)=0\}$.  In the regular part $\Sg\setminus\Sg_0$,
we can define the \emph{horizontal unit normal vector} $\nu_H$ by
$\nu_H:=N_H/|N_H|$.  Consider the unit vector field $Z$ on $\Sg\setminus\Sg_0$
given by $Z:=J(\nu_H)$.  As $Z$ is horizontal and orthogonal to
$\nu_H$, it follows that $Z$ is tangent to $\Sg$.  The integral curves
of $Z$ in $\Sg\setminus\Sg_{0}$ will be called \emph{characteristic curves}.
Characteristic curves foliate the regular part of $\Sg$.

%\subsection{Mean curvature. Variation formulas}

Consider a $C^1$ vector field $U$ with compact support on $\hhn$, and
denote by $\{\varphi_t\}_{t\in\rr}$ the associated group of
diffeomorphisms.  Let $E$ be a bounded region enclosed by a
hypersurface $\Sg$.  The families $\{E_t\}$, $\{\Sg_{t}\}$, for $t$
small, are the variations of $E$ and $\Sg$ induced by $U$.  Let
$V(t)=|E_t|$ and $A(t)=A(\Sg_t)$.  We say that the variation is
\emph{volume-preserving} if $V'(0)=0$.  We say that $\Sg$ is
\emph{area-stationary} if $A'(0)=0$ for any variation, and
\emph{volume-preserving area-stationary} if $A'(0)=0$ for any volume
preserving variation.

If $\Sg$ is a $C^1$ hypersurface enclosing a bounded region $E$, it 
is well-known that
\begin{equation}
\label{eq:1stvol}
V'(0)=\int_E\divv U\,dv=-\int_\Sg u\,d\Sg,
\end{equation}
where $u=\escpr{U,N}$ and $N$ is the unit vector normal to $\Sg$
pointing into $E$.

If $\Sg$ is $C^2$, and $N$ is a unit vector field normal to $\Sg$, the
\emph{mean curvature} of $\Sg\setminus\Sg_{0}$ is given by
$-nH:=\divv_{\Sg}\nuh$, where
$\divv_{\Sg}U(p):=\sum_{i=1}^{2n}\escpr{D_{e_{i}}U,e_{i}}$ for any
orthonormal basis $\{e_{i}\}$ of $T_{p}\Sg$.  We say that $\Sg$ has
constant mean curvature if $H$ is constant on $\Sg\setminus\Sg_{0}$.  By
combining \cite[Lemma~3.2]{rr1} and \cite[Lemma~4.3 (4.7)]{rr2} we have

\begin{lemma}
\label{lem:1starea}
Let $\Sg\subset\hh^n$ be a $C^2$ hypersurface enclosing a bounded
region $E$, with inner unit normal vector $N$.  Consider a variation
induced by a vector field $U$, and let $u=\escpr{U,N}$.  Assume that
$\Sg$ is volume-preserving area-stationary, and let $H$ be the mean
curvature of $\Sg$.  Then we have
\[
A'(0)=-\int_{\Sg}nH\,u\,d\Sg.
\]
\end{lemma}
We remark that the variation associated to $U$ in the statement of
Lemma~\ref{lem:1starea} is not assumed to be volume-preserving.  From
Lemma~\ref{lem:1starea} and \eqref{eq:1stvol} we easily obtain
\begin{lemma}
\label{lem:a-nhv}
Let $\Sg\subset\hh^n$ be a $C^2$ hypersurface enclosing a bounded
region $E$.  Assume that $\Sg$ is volume-preserving area-stationary,
and let $H$ be the $($constant$)$ mean curvature of $\Sg$.  Then $\Sg$ is a critical
point of the functional $A-nHV$ for any variation.
\end{lemma}

Let $\Sg$ be a two-sided $C^2$ hypersurface without singular points.
We can translate it vertically to get a foliation of the vertical
cylinder $C$ over $\Sg$.  Denote by $N$ the unit normal to the
foliation, and by $\nuh$ the horizontal unit normal obtained from $N$.
For any $p\in C$, let $\{e_{i}\}$ be an orthonormal basis of the leaf
passing through $p$.  Then
\[
\divv\nuh(p)=\sum_{i=1}^{2n}\escpr{D_{e_{i}}\nuh,e_{i}}
+\escpr{D_{N_{p}}\nuh,N_{p}}
=-nH(p)+\escpr{D_{N_{p}}\nuh,N_{p}}.
\]
Since $N=|N_{H}|\,\nuh+\escpr{N,T}\,T$, and 
$\escpr{D_{U}\nuh,\nuh}$, $\escpr{D_{T}\nuh,T}$, and
$\escpr{D_{\nuh}T,\nuh}=\escpr{J(\nuh),\nuh}$ vanish, we conclude 
that $\escpr{D_{N_{p}}\nuh,N_{p}}=0$. Hence
\begin{equation}
\label{eq:divv}
-nH=\divv\nuh.
\end{equation}

In \cite[Lemma~3.1 (3.3)]{rr1} it is proven that
\begin{equation}
\label{eq:dunuh}
D_{u}\nu_H=|N_H|^{-1}\,\sum_{i=1}^{2n-1}\big(\escpr{D_uN,z_i}-\escpr{N,T}\,
\escpr{J(u),z_i}\big)\,z_i+\escpr{z,u}\,T,
\end{equation}
where $u\in T_{p}\Sg$ for a given point $p\in\Sg$, and
$\{z_{1},\ldots,z_{2n-1}\}$ is an orthonormal basis of
$T_{p}\Sg\cap\mathcal{H}_{p}$, with $z_{1}=z=J((\nuh)_{p})$.
Completing $\{z_{i}\}$ to an orthonormal basis of $T_{p}\Sg$ by adding
a vector $v$, we obtain from \eqref{eq:dunuh} that
$\escpr{D_{v}\nuh,v}=0$.  Hence we conclude
\[
-nH(p)=\sum_{i=1}^{2n-1}\escpr{D_{z_{i}}\nuh,z_{i}},
\]
where $\{z_{i}\}$ is an orthonormal basis of
$T_{p}\Sg\cap\mathcal{H}_{p}$.  From \eqref{eq:dunuh}, it follows that
the endomorphism $v\mapsto -D_{v}\nuh-|N_H|^{-1}\escpr{N,T}\,J(v)^\top$, defined in the
subspace $T_{p}\Sg\cap\mathcal{H}_{p}$, $p\in\Sg\setminus\Sg_{0}$, is
selfadjoint.  Thus there exists an orthonormal basis
$\{v_{1},\ldots,v_{2n-1}\}$ of $T_{p}\Sg\cap\mathcal{H}_{p}$ composed
of eigenvectors with eigenvalues $\kappa_{1},\ldots,\kappa_{2n-1}$.
By analogy with the Riemannian case they will be named \emph{principal
curvatures}, and we have $-nH=\kappa_{1}+\ldots+\kappa_{2n-1}$.

\subsection{Geodesics in $\hh^n$}

We refer the reader to \cite[\S~3]{rr1} for detailed arguments.
Geo\-desics in $\hh^n$ are horizontal curves $\ga:I\to\hh^n$ which are
critical points of the Riemannian length $L(\ga):=\int_{I}|\dot{\ga}|$
for any variation by horizontal curves $\ga_{\eps}$.  A vector field
$U$ along $\ga$ induces a variation by horizontal curves if and only
if
\begin{equation}
\label{eq:admissible}
\dot{\ga}\big(\escpr{U,T}\big)
+2\,\escpr{\dot{\ga},J(U)}=0.
\end{equation}
The derivative of length for such a variation is given by
\begin{equation}
\label{eq:1stlength}
\frac{d}{d\eps}\bigg|_{\eps=0} 
L(\ga_{\eps})=-\int_{I}\escpr{D_{\dot{\ga}}\dot{\ga},U}.
\end{equation}
Observe that $D_{\dot{\ga}}\dot{\ga}$ is orthogonal to both
$\dot{\ga}$ and $T$.  Along $\ga$ consider the orthonormal basis of
$T\hh^n$ given by $T$, $\dot{\ga}$, $J(\dot{\ga})$, $Z_{1}$, $\ldots$,
$Z_{2n-2}$.  In the same way as for the case of $\hh^1$, see
\cite[\S~3]{rr2}, we take any smooth $f:I\to\rr$ vanishing at the
endpoints of $I$ so that $\int_{I}f=0$.  The vector field $U$ along
$\ga$ so that $U_{H}=f\,J(\dot{\ga})$, and $\escpr{U,T}=2\,\int_{I}
f$, satisfies \eqref{eq:admissible}.  Hence \eqref{eq:1stlength} allows
us to conclude that $\escpr{D_{\dot{\ga}}\dot{\ga},J(\dot{\ga})}$ is
constant.  Now let $f:I\to\rr$ be any smooth function vanishing at the
endpoints of $I$.  Then the vector field $U=f\,Z_{i}$, for any
$i=1,\ldots,2n-2$, satisfies \eqref{eq:admissible}, and hence
$D_{\dot{\ga}}\dot{\ga}$ is orthogonal to $Z_{i}$ for all
$i=1,\ldots,2n-2$.  So we obtain that the horizontal geodesic
$\ga:I\to\hh^n$ satisfies the equation
\begin{equation}
\label{eq:geodesic}
D_{\dot{\ga}}\dot{\ga}+2\la\,J(\dot{\ga})=0,
\end{equation}
for some constant $\la\in\rr$. For $\la\in\rr$, $p\in\hh^n$, and 
$v\in T_{p}\hh^n$, $|v|=1$, the geodesic $\ga:I\to\hh^n$ of curvature 
$\la$ with initial conditions $\ga(0)=p$, $\dot{\ga}(0)=v$, will be 
denoted by $\ga_{p,v}^\la$.

The equations of a geodesic can be computed in coordinates in the 
following way: let 
$\ga(s)=(x_{1}(s),y_{1}(s),\ldots,x_{n}(s),y_{n}(s),t(s))$ be a 
horizontal geodesic. Then
\[
\dot{\ga}(s)=\sum_{i=1}^n \dot{x}_{i}(s)\,(X_{i})_{\ga(s)}
+\dot{y}_{i}(s)\,(Y_{i})_{\ga(s)},
\]
and
\[
\dot{t}(s)=\sum_{i=1}^n (\dot{x}_{i} y_{i}-x_{i}\dot{y}_{i})(s).
\]
So equation \eqref{eq:geodesic} is transformed in
\begin{align*}
\ddot{x}_{i}&=2\la\,\dot{y}_{i},
\\
\ddot{y}_{i}&=-2\la\,\dot{x}_{i},
\end{align*}
with initial conditions $x_{i}(0)=(x_{0})_{i}$,
$y_{i}(0)=(y_{0})_{i}$, and $\dot{x}_{i}(0)=A_{i}$,
$\dot{y}_{i}(0)=B_{i}$, with $\sum_{i=1}^n (A_{i}^2+B_{i}^2)=1$.

Integrating these equations, for $\la=0$, we obtain
\begin{align*}
x_{i}(s)&=(x_{0})_{i}+A_{i}s, \\
y_{i}(s)&=(y_{0})_{i}+B_{i}s, \\
t(s)&=t_{0}+\sum_{i=1}^n (A_{i}(y_{0})_{i}-B_{i}(x_{0})_{i})\,s,
\end{align*}
which are horizontal Euclidean straight lines in $\hh^n$.

Integrating, for $\la\neq 0$, we obtain
\begin{align*}
x_{i}(s)&=(x_{0})_{i}+A_{i}\,\bigg(\frac{\sin(2\la s)}{2\la}\bigg)
+B_{i}\,\bigg(\frac{1-\cos(2\la s)}{2\la}\bigg),
\\
y_{i}(s)&=(y_{0})_{i}-A_{i}\,\bigg(\frac{1-\cos(2\la s)}{2\la}\bigg)
+B_{i}\,\bigg(\frac{\sin(2\la s)}{2\la}\bigg),
\\
t(s)&=t_{0}+\frac{1}{2\la}\,\bigg(s-\frac{\sin(2\la s)}{2\la}\bigg)
+\sum_{i=1}^n\bigg\{
(A_{i}(x_{0})_{i}+B_{i}(y_{0})_{i})\bigg(\frac{1-\cos(2\la s)}{2\la}\bigg)
\\
&\qquad-(B_{i}(x_{0})_{i}-A_{i}(y_{0})_{i})
\bigg(\frac{\sin(2\la s)}{2\la}\bigg)\bigg\}.
\end{align*}
In case $x_{0}=y_{0}=0$, we obtain
\begin{align*}
x_{i}(s)&=A_{i}\,\bigg(\frac{\sin(2\la s)}{2\la}\bigg)
+B_{i}\,\bigg(\frac{1-\cos(2\la s)}{2\la}\bigg),
\\
y_{i}(s)&=-A_{i}\,\bigg(\frac{1-\cos(2\la s)}{2\la}\bigg)
+B_{i}\,\bigg(\frac{\sin(2\la s)}{2\la}\bigg),
\end{align*}
and so $\dot{x}_{i}$, $\dot{y}_{i}$, $i=1,\ldots,n$, can be expressed
in terms of $x_{i}$, $y_{i}$ in the following way
\begin{align}
\label{eq:x0}
\dot{x}_{i}(s)&=\frac{\la\,\sin(2\la s)}{1-\cos(2\la s)}\,x_{i}(s)
+\la\,y_{i}(s), \\
\label{eq:y0}
\dot{y}_{i}(s)&=-\la\,x_{i}(s)
+\frac{\la\,\sin(2\la s)}{1-\cos(2\la s)}\,y_{i}(s).
\end{align}

\subsection{The spheres $\sph_{\la}$}

For any $\la>0$, $p\in\hh^n$, consider the hypersurface 
$\sph_{\la,p}$ defined by
\[
\sph_{\la,p}:=\bigcup_{v\in\mathcal{H}_{p}, 
|v|=1}\ga_{p,\la}^v([0,\pi/\la]).
\]
If $p$ is translated to the point $q=(0,-\pi/(4\la^2))$, then 
$\sph_{\la}:=\sph_{\la,q}$ is the union of the graphs associated to 
the functions $f$ and $-f$, where
\[
f(z)=\frac{1}{2\la^2}\,\{\la 
|z|\,\sqrt{1-\la^2\,|z|^2}+\arccos{(\la |z|)}\},
\qquad |z|\leq\frac{1}{\la}.
\]
The hypersurface $\sph_\la$ is compact and homeomorphic to a
$(2n)$-dimensional sphere.  Its singular set consists of the two points
$\pm(0,\pi/(4\la^2))$ on the $t$-axis, called the \emph{poles}.  It
is known that the spheres $\sph_\la$ are $C^2$ but not $C^3$ around
the singular points.  These hypersurfaces were conjectured to be the
(smooth) solutions to the isoperimetric problem in $\hh^1$ by P.~Pansu
\cite{pansu1}.  It was proven in \cite{rr1} that the
hypersurfaces $\sph_{\la}$ are the only compact hypersurfaces of
revolution with constant mean curvature $\la$ in $\hh^n$.  We shall
denote by $\bb_\la$ the topological closed ball enclosed by
$\sph_\la$.  It is well known that the spheres $\sph_{\la}$ consist of
the union of segments of geodesics of curvature $\la$ and length
$\pi/\la$ starting from a given point $p\in\hh^n$.

\begin{lemma}
\label{lem:char}
The characteristic curves in $\sph_{\la}$ are the geodesics of 
curvature $\la$ joining the poles.
\end{lemma}

\begin{proof}
Since $\sph_{\la}^\pm$ is the graph of the function $\pm f$, the inner
unit normal to $\sph_{\la}$ is proportional to
\[
\sum_{i=1}^n\bigg\{\bigg(\frac{\ptl f}{\ptl x_{i}}-y_{i}\bigg)\,X_{i}
+\bigg(\frac{\ptl f}{\ptl y_{i}}+x_{i}\bigg)\,Y_{i}\bigg\}-T
\]
on $\sph_{\la}^-$, and proportional to
\[
-\sum_{i=1}^n\bigg\{\bigg(\frac{\ptl (-f)}{\ptl x_{i}}-y_{i}\bigg)\,X_{i}
+\bigg(\frac{\ptl (-f)}{\ptl y_{i}}+x_{i}\bigg)\,Y_{i}\bigg\}+T
\]
on $\sph_{\la}^-$.  Let $\nu_{\la}$ be the horizontal unit normal to
$\sph_{\la}$.  From
\[
\frac{\ptl f}{\ptl x_{i}}=-\la |z|\,(1-\la^2|z|^2)^{-1/2}\,x_{i},
\qquad
\frac{\ptl f}{\ptl y_{i}}=-\la |z|\,(1-\la^2|z|^2)^{-1/2}\,y_{i},
\]
we have that $J(\nu_{\la})$ is given by
\begin{equation*}
\sum_{i=1}^n\bigg(\!-\la x_{i}
-\frac{(1-\la^2|z|^2)^{1/2}}{|z|}\,y_{i}\bigg)\,Y_{i}
+\bigg(\la y_{i}-\frac{(1-\la^2|z|^2)^{1/2}}{|z|}\,x_{i}\bigg)\,X_{i}
\end{equation*}
on $\sph_{\la}^+$, and by
\[
\sum_{i=1}^n \bigg(\!-\la x_{i}
+\frac{(1-\la^2|z|^2)^{1/2}}{|z|}\,y_{i}\bigg)\,Y_{i}
+\bigg(\la y_{i}+\frac{(1-\la^2|z|^2)^{1/2}}{|z|}\,x_{i}\bigg)\,X_{i}
\]
on $\sph_{\la}^-$, where $|z|^2=\sum_{i=1}^n (x_{i}^2+y_{i}^2)$.

On the other hand, the tangent vector to a horizontal geodesic of
curvature $\la$ leaving from $(0,-\pi/(4\la^2))$, is given by
$\dot{\ga}(s)=\sum_{i=1}^n
\dot{x}_{i}(s)\,X_{i}+\dot{y}_{i}(s)\,Y_{i}$, where $\dot{x}_{i}(s)$,
$\dot{y}_{i}(s)$ satisfy \eqref{eq:x0} and \eqref{eq:y0}.  A direct
computation shows
\begin{equation}
\label{eq:la}
\frac{\la\,\sin(2\la s)}{1-\cos(2\la s)}=\pm\frac{(1-\la^2|z|^2)^{1/2}}{|z|},
\end{equation}
where the plus sign is chosen in case $s\in [0,\pi/(2\la)]$, that is,
when $\ga(s)\in\sph_{\la}^-$, and the minus sign if $s\in
[\pi/(2\la),\pi/\la]$, when $\ga(s)\in\sph_{\la}^+$.  Replacing the
value of $\la\sin(2\la s)/(1-\cos(2\la s))$ in equations
\eqref{eq:x0} and \eqref{eq:y0} by using \eqref{eq:la}, we conclude 
that $\dot{\ga}$ is equal to $J(\nu_{\la})$, i.e., $\ga$ is a 
characteristic curve of $\sph_{\la}$.
\end{proof}

\begin{lemma}
\label{lem:spheres}
Let $\Sg\subset\hh^n$ be a $C^2$ compact hypersurface with a finite
number of isolated singular points.  Assume that $\Sg$ has constant
mean curvature.  Then $\Sg$ is volume-preserving area-stationary.  In
particular, the spheres $\sph_{\la}$ are volume-preserving
area-stationary.
\end{lemma}

\begin{proof}
Let $U$ be a vector field inducing a volume-preserving variation of
$\Sg$.  Let $u=\escpr{U,N}$.  By the first variation of volume
\eqref{eq:1stvol} we have $\int_{\Sg} u\,d\Sg=0$.  By the first 
variation of area \cite[Lemma~3.2]{rr1}, we have
\[
A'(0)= -\int_\Sg\divv_\Sg\big(u\,(\nu_H)^\top\big)\,d\Sg,
\]
since $u$ has mean zero and $\divv_{\Sg}\nuh$ is constant.

To analyze the above integral, we consider open balls
$B_{\eps}(p_{i})$ of radius $\eps>0$ centered at the points $p_{1}$,
$\ldots$, $p_{k}$ of the singular set $\Sg_{0}$.  By the divergence
theorem in $\Sg$, we have, for $\Sg_{\eps}=\Sg\setminus\bigcup_{i=1}^k
B_{\eps}(p_{i})$,
\begin{equation*}
-\int_{\Sg_{\eps}}\divv_\Sg\big(u\,(\nu_H)^\top\big)\,d\Sg=
\sum_{i=1}^k\,
\int_{\ptl B_{\eps}(p_{i})}
u\,\escpr{\xi_{i},(\nuh)^\top}\,d(\ptl B_{\eps}(p_{i})),
\end{equation*}
where $\xi_{i}$ is the inner unit normal vector to $\ptl
B_{\eps}(p_{i})$ in $\Sg$, and $d(\ptl B_{\eps}(p_{i}))$ is the
Riemannian volume element of $\ptl B_{\eps}(p_{i})$.  Note also that
\[
\bigg|\sum_{i=1}^k\,\int_{\ptl
B_{\eps}(p_{i})}u\,\escpr{\xi_{i},(\nu_{H})^\top}\,d(\ptl B_{\eps}(p_{i}))
\bigg|
\leq\big(\sup_{\Sg}\,|u|\big)\,\sum_{i=1}^k\vol_{2n-1}(\ptl B_{\eps}(p_{i})),
\]
where $\vol_{2n-1}(\ptl B_{\eps}(p_{i}))$ is the Riemannian
$(2n-1)$-volume of $\ptl B_{\eps}(p_{i})$.  Observe that the function
$|\divv_{\Sg}(u\,(\nu_{H}^\top))|$ is bounded from above by
$(\sup_{\Sg}\,|u|)\,|\divv_{\Sg}\nu_{H}
-|N_{H}|\,|\divv_{\Sg}N|+|\nabla_{\Sg}u|$, which is also bounded.  So
we can apply the dominated convergence theorem and the fact that
$\vol_{2n-1}(\ptl B_{\eps}(p_{i}))\to 0$, when $\eps\to 0$, to prove
that $A'(0)=0$.
\end{proof}

% \subsection{Finite perimeter sets in $\hh^n$}

\subsection{$\hh$-regular surfaces (\cite{ascv}, \cite{fssc})}

Let $\Om\subset\hh^n$ be an open set.  Then $C^1_{\hh}(\Om)$ is the
set of continuous real functions in $\Om$ such that $\nabla_{\hh}f$ is
continuous, \cite[def.~5.7]{fssc}, where $\nabla_{\hh}f$ is defined by
\[
\nabla_{\hh}f:=\sum_{i=1}^n X_{i}(f)\,X_{i}+Y_{i}(f)\,Y_{i}.
\]
Following \cite[def.~6.1]{fssc} we say that $\Sg\subset\hh^n$ is an
$\hh$-regular hypersurface if, for every $p\in\Sg$, there is an open
set $\Om$ containing $p$, and a function $f\in C^1_{\hh}(\Om)$ such
that
\[
\Sg\cap\Om=\{q\in\Om : f(q)=0\}, \qquad \nabla_{\hh}f(q)\neq 0.
\]
From \cite[Thm.~6.5]{fssc}, we know that if $\Sg$ is an $\hh$-regular 
hypersurface defined locally by a function $f\in C^1_{\hh}(\Om)$ with 
nonvanishing horizontal gradient $\nabla_{\hh}f$, and we let
$E:=\{q\in\Om : f(q)<0\}$, then $E$ is a set of locally finite perimeter in
$\Om$ and $\nuh=-\nabla_{\hh}f/|\nabla_{\hh}f|$.

We have the following

\begin{lemma}
\label{lem:charcurve}
Let $\Sg$ be an $\hh$-regular hypersurface, and let $\nuh$ be its horizontal unit~normal vector. If $\ga:I\to\hh^n$ is an integral curve of $J(\nu_H)$ with $\ga(0)\in\Sg$, then $\ga(I)\subset\Sg$.
\end{lemma}

\begin{proof}
Locally $\Sg=\{p\in\Om : f(p)=0\}$, for some open set
$\Om\subset\hh^n$, and some $f\in C^1_{\hh}(\Om)$ with non vanishing
horizontal gradient. It is enough to prove the result in $\Sg\cap\Om$. If $\ga$ is any $C^1$ curve then
\begin{equation}
\label{eq:chain}
\frac{d}{ds} f(\ga(s))=\escpr{\nabla_{\hh}f(\ga(s)),\dot{\ga}(s)}.
\end{equation}
Let us check this equality. The $C_\hh^1$ function $f$ is Pansu-di\-ffe\-rentiable in its domain of definition by \cite[Thm.~5.10]{fssc}. By the results in \cite[\S~5]{fssc} this implies that, for $p$ in a fixed small neighborhood of $a$, we have
\begin{equation}
\label{eq:taylor}
f(p)-f(a)-\escpr{\nabla_\hh f(a),\pi_a(a^{-1}\cdot p))}=||a^{-1}\cdot p||_\infty E(a,p),
\end{equation}
where the function $E(a,p)$ satisfies
\begin{equation}
\label{eq:limE}
\lim_{p\to a} E(a,p)=0,
\end{equation}
$||\cdot||_\infty$ is the homogeneous norm defined in \eqref{eq:norm}, and the projection $\pi_a$ is defined, as in \cite[Def.~2.19]{fssc}, by
\[
\pi_a(z,t):=\sum_{i=1}^n x_i\,(X_i)_a+y_i\,(Y_i)_a.
\]
Hence, taking $a=\ga(s)$, $p=\ga(s+h)$, $g(s)=(f\circ\ga)(s)$, we get from \eqref{eq:taylor}
%\begin{equation}
\label{eq:taylor2}
\begin{multline}\frac{g(s+h)-g(s)}{h}=\escpr{\nabla_\hh f(\ga(s)),\frac{\pi_{\ga(s)}(\ga(s)^{-1}\cdot\ga(s+h))}{h}}\\
+
\frac{||\ga(s)^{-1}\cdot\ga(s+h)||_\infty}{h}\,E(\ga(s),\ga(s+h)).
\end{multline}
%\end{equation}
Since $\ga$ is a horizontal curve, it is immediate to check that
\begin{equation}
\label{eq:lim1}
\lim_{h\to 0}\frac{\pi_{\ga(s)}(\ga(s)^{-1}\cdot\ga(s+h))}{h}=\dot{\ga}(s).
\end{equation}
On the other hand, since the norm $||\cdot||_\infty$ is globally equivalent to the Carnot-Ca\-ra\-th\'eodory distance, there is a constant $C>0$ so that $||\ga(s)^{-1}\cdot\ga(s+h)||_\infty\le C\,d(\ga(s),\ga(s+h))$. As $\ga$ is a $C^1$ horizontal curve we have $d(\ga(s),\ga(s+h))\le ||\dot{\ga}||_\infty\,h$. So we have that $||a^{-1}\cdot\ga(s+h)||_\infty/h$ is bounded and
\begin{equation}
\label{eq:lim2}
\lim_{h\to 0}\frac{||\ga(s)^{-1}\cdot\ga(s+h)||_\infty}{h}\,E(\ga(s),\ga(s+h))=0,
\end{equation}
by \eqref{eq:limE}. Taking limits when $h$ goes to $0$ in \eqref{eq:taylor2}, and using \eqref{eq:lim1} and \eqref{eq:lim2}, we get \eqref{eq:chain}.

Assume now that $\ga$ is an integral curve of $J(-\nabla_\hh f/|\nabla_\hh f|)$ with $\ga(0)\in\Sg$. Since $\nabla_\hh f$ and $J(\nabla_\hh f)$ are orthogonal, equation \eqref{eq:chain} implies
\[
\frac{d}{ds} f(\ga(s))=0.
\]
We conclude that $f\circ\ga$ is a constant function and, since $f(\ga(0))=0$, we have $f\circ\ga\equiv 0$ in $\Om$, and so $\ga(I)\subset\Sg$. Since $\nu_H=-\nabla_\hh f/|\nabla_\hh f|$ on $\Sg$, $\ga$ in an integral curve of $J(\nu_H)$ contained in $\Sg$.
\end{proof}

% \begin{lemma}
% \label{lem:hregular}
% Let $\Om\subset\hh^n$ be an open set, and let $E\subset\Om$ be a
% finite perimeter set. Assume that, for every $p\in\ptl E\cap\Om$,
% $\nuh(p)$ exists and it is continuous. Then $\ptl E\cap\Om$ is an
% $\hh$-regular hypersurface.
% \end{lemma}
% 
% \begin{proof}
% Following \cite{fssc}.
% \end{proof}

\section{Proof of the isoperimetric inequality}
\label{sec:isop}

We shall denote by $D_r:=\{(z,0) :|z|\le r\}$ the closed Euclidean
disk of radius $r>0$ contained in the Euclidean hyperplane
$\Pi_0:=\{t=0\}$, and by $C_r:=\{(z,t) : |z|\le r\}$ the vertical
cylinder over $D_r$.  The vertical $t$-axis $\{(0,t): t\in\rr\}$
will be denoted by $L$.  For any set $B\subset\hh^n$, define
$B^+:=B\cap\{(z,t)\in\hh^n : t\ge 0\}$, $B^-:=B\cap\{(z,t)\in\hh^n:
t\le 0\}$.

We recall that to any finite perimeter set $E\subset\hh^n$ we may add its density one points and remove its density zero points, without changing the perimeter and the volume of $E$. We shall always normalize a finite perimeter set in this way.

\begin{theorem}
\label{th:main}
Let $E\subset\hh^n$ be a finite perimeter set such that $D\subset E\subset C$, where $D=D_r$, $C=C_r$, for some $r>0$. Then
\begin{equation}
\label{eq:main}
|\ptl E|\ge |\ptl\bb_{\mu}|,
\end{equation}
where $\bb_{\mu}$ is the ball with $|\bb_{\mu}|=|E|$.  Equality holds in \eqref{eq:main} if and only if $\bb_{\mu}=\overline{E}$.
\end{theorem}

%\begin{theorem}
%\label{th:main}
%For some $r>0$, let $C=C_{r}$, $D=D_{r}$. Let $E\subset\hh^n$ be a
%normalized finite perimeter set such that $D\subset E\subset C$. Then
%\begin{equation}
%\label{eq:main}
%|\ptl E|\ge |\ptl\bb_{\mu}|,
%\end{equation}
%where $\bb_{\mu}$ is the ball with $|\bb_{\mu}|=|E|$.  In case $\ptl
%E$ is locally lipschitz in Euclidean sense, or $\ptl E$ is an
%$\hh$-regular hypersurface, then equality holds in \eqref{eq:main} if
%and only if $E=\bb_{\mu}$.
%\end{theorem}

\begin{proof}
It can be easily proven that $E^\pm:=E\cap (\hh^n)^\pm$ are finite perimeter sets.  The reduced boundary $\ptl^*E^+$ of $E^+$ is contained in $(\ptl^*E\cap\{t>0\})\cup\intt(D)$, where $\intt(D)$ is the interior of $D$ inside $\Pi_{0}$.

We choose two families of functions. For $0<\eps<1$ we consider smooth functions $\varphi_{\eps}$,~depending on the Euclidean distance to the vertical axis $L$, so that $0\le\varphi_\eps\le 1$, and
\begin{align*}
\varphi_{\eps}(p)= 0, \qquad & d(p,L)\le\eps^2, \\
\varphi_{\eps}(p)= 1, \qquad & d(p,L)\ge\eps, \\
|\nabla\varphi_{\eps}(p)|\le 2/\eps, \qquad& \eps^2\le d(p,L)\le\eps.
\end{align*}
Again for $0<\eps<1$ we consider smooth functions $\psi_{\eps}$, depending on
the distance to the Euclidean hyperplane $\Pi_{0}$, so
that $0\le\psi_\eps\le 1$, and
\begin{align*} 
\psi_{\eps}(p)=1, \qquad & d(p,\Pi_{0})\le \eps^{-1/2}, \\
\psi_{\eps}(p)=0, \qquad & d(p,\Pi_{0})\ge \eps^{-1/2}+1, \\
|\nabla\psi_{\eps}(p)|\le 2, \qquad&\eps^{-1/2}\le d(p,\Pi_{0})\le
\eps^{-1/2}+1.
\end{align*}

Let $\la:=1/r$.  Then the ball $\bb_{\la}$ satisfies
$\bb_{\la}\cap\Pi_0=D$.  Translate vertically the closed halfspheres
$\sph_{\la}^+$ to get a foliation of $C$.  Let $X$ be the vector field
on $C\setminus L$ given by the horizontal unit normal to the leaves of the
foliation.  By \eqref{eq:divv}, on $C\setminus L$ we have
\[
\divv X=-n\la.
\]

We consider the horizontal vector field $\psi_{\eps}\varphi_{\eps}X$,
which has compact support on $\hh^n$.  
\[
\int_{E^+}\divv(\psi_{\eps}\varphi_{\eps}X)\,dv=
\int_{E^+}\psi_{\eps}\varphi_{\eps}\divv X\,dv+
\int_{E^+}\escpr{\nabla(\psi_{\eps}\varphi_{\eps}),X}\,dv.
\]
Observe that
\[
\lim_{\eps\to 0} \int_{E^+}\psi_{\eps}\varphi_{\eps}\divv 
X\,dv=-n\la\,|E^+|,
\]
by Lebesgue's Dominated Convergence Theorem since
$\psi_{\eps}\varphi_{\eps}\divv X$ is uniformly bounded, $E^+$ has
finite volume, and $\lim_{\eps\to 0} \psi_\eps\varphi_\eps =1$. On the other hand
\[
\lim_{\eps\to
0}\int_{E^+}\escpr{\nabla(\psi_{\eps}\varphi_{\eps}),X}\,dv=0,
\]
since $\escpr{\varphi_{\eps}\nabla\psi_{\eps},X}$ is bounded and
converges pointwise to $0$, and
\[
\lim_{\eps\to 0}\int_{E^+}|\escpr{\psi_{\eps}\nabla\varphi_{\eps},X}| 
\,dv\le 
\lim_{\eps\to 0}\int_{(\hh^n)^+}\psi_\eps\,|\nabla\varphi_{\eps}|\,dv=0.
\]
The last equality is easily checked taking classical cilindrical coordinates in $\hh^n=\rr^{2n+1}$. So we conclude
\begin{equation}
\label{eq:vol+}
\lim_{\eps\to 0}\int_{E^+}
\divv(\psi_{\eps}\varphi_{\eps}X)\,dv=-n\la\,|E^+|.
\end{equation}

By applying the Divergence Theorem for finite perimeter sets \cite{fssc} to $E^+$ and to the
vector field $\psi_{\eps}\varphi_{\eps}X$ we have
\[
\int_{E^+}\divv(\psi_{\eps}\varphi_{\eps}X)\,dv=
-\int_{D}\escpr{\psi_{\eps}\varphi_{\eps}X, N_{D}}\,dD
-\int_{\ptl^*E\cap\{t>0\}}\escpr{\psi_{\eps}\varphi_{\eps}X,\nuh}
\,d\,|\ptl E|,
\]
where $N_{D}$ is the Riemannian unit normal to $D$ pointing into $(\hh^n)^+$, $dD$ is the Riemannian area element on $D$, and $\nuh$ is the inner horizontal unit normal to $\ptl^*E$. Taking limits when $\eps\to 0$, we get from \eqref{eq:vol+} and inequality $\escpr{X,\nuh}\le 1$,
\begin{equation}
\label{eq:om+}
-n\la\,|E^+|\ge
-\int_{D}\escpr{X,N_{D}}\,dD-\int_{\ptl^*E\cap\{t> 0\}} d\,|\ptl E|,
\end{equation}
with equality if and only if, $|\ptl E|$-a.e., $X=\nuh$ on $\ptl E\cap (\hh^n)^+$.

We may replace $E^+$ by $\bb_{\la}^+$ in the previous reasoning 
to obtain
\begin{equation}
\label{eq:bbla+}
-n\la\,|\ptl\bb_{\la}^+|=
-\int_{D}\escpr{X,N_{D}}\,dD-\int_{\ptl^*\bb_{\la}\cap\{t>0\}}
d\,|\ptl\bb_{\la}|.
\end{equation}
Hence from \eqref{eq:om+} and \eqref{eq:bbla+} we obtain
\begin{equation}
\label{eq:+}
\int_{\ptl^*E\cap\{t>0\}}d\,|\ptl E|\ge
\int_{\ptl\bb_{\la}\cap\{t>0\}} d\,|\ptl\bb_{\la}|
+n\la\,\big(|E^+|-|\bb_{\la}^+|\big),
\end{equation}
with equality if and only if, $|\ptl E|$-a.e., $X=\nuh$ on $\ptl E\cap (\hh^n)^+$.

We consider now the foliation of $C$ by vertical translations of the
closed halfspheres $\sph_\la^-$.  Let $Y$ be the vector field on $C\setminus L$
given by the horizontal unit normal to the leaves of the foliation.
By applying the previous argument we get a similar estimate
\begin{equation}
\label{eq:-}
\int_{\ptl^*E\cap\{t<0\}}d\,|\ptl E|\ge
\int_{\ptl\bb_{\la}\cap\{t<0\}} d\,|\ptl\bb_{\la}|
+n\la\,\big(|E^-|-|\bb_{\la}^-|\big),
\end{equation}
with equality if and only if, $|\ptl E|$-a.e., $Y=\nuh$ on 
$\ptl E\cap (\hh^n)^-$.

Hence, adding \eqref{eq:+} and \eqref{eq:-}, and taking into account $|\ptl E|\ge |\ptl E|(\hh^n\setminus\Pi_0)$, that $\ptl\bb_{\la}\cap\Pi_0$ do not contribute to the perimeter of $E$ and $\bb_\la$, and that $E\cap\Pi_{0}$ and $\bb_{\la}\cap\Pi_{0}$ do not contribute to the volume of $E$ and $\bb_{\la}$, we get
\begin{equation}
\label{eq:estimate}
|\ptl E|\ge  |\ptl\bb_{\la}|+n\la\,\big(|E|-|\bb_{\la}|\big),
\end{equation}
with equality if and only if, $|\ptl E|$-a.e, $X=\nuh$ on $\ptl E\cap (\hh^n)^+$ and
$Y=\nuh$ on $\ptl E\cap (\hh^n)^-$ and $|\ptl E|=|\ptl E|(\hh^n\setminus\Pi_0)$.

Let $f(\rho):=n\rho\,|E|+ |\ptl\bb_\rho|-n\rho\,|\bb_\rho|$.
By Lemmae~\ref{lem:a-nhv} and \ref{lem:spheres}, the sphere $\sph_\rho$
is a critical point of $A-n\rho V$, with $\rho$ fixed, for any
variation.  So we have $\area(\sph_\rho)'-n\rho\,\vol(\bb_\rho)'=0$,
where primes indicate the derivative with respect to $\rho$.  Hence we
have
\[
f'(\rho)=n\,(|E|-|\bb_\rho|).
\]
Since the function $\rho\mapsto |\bb_\rho|$ is strictly decreasing and takes
its values in the interval $(0,+\infty)$, we obtain that $f(\rho)$ is an
strictly convex function with a unique minimum $\mu$ for which
$|E|=|\bb_{\mu}|$.  Hence we obtain from \eqref{eq:estimate}
\begin{equation}
\label{eq:main2}
|\ptl E|\ge f(\la)\ge f(\mu)=|\ptl\bb_{\mu}|,
\end{equation}
which implies \eqref{eq:main}.

Assume now that equality holds in \eqref{eq:main2}. Then, since $f$ is strictly convex, $\la=\mu$.~By \cite[Thm~1.2]{mv}, $\ptl E\setminus L$ is an $\hh$-regular hypersurface. By Lemma~\ref{lem:charcurve}, the integral curves of $J(\nu_{H})$, starting from points in $\ptl E\setminus L$, are contained in $\ptl E\setminus L$.

Observe now that $\ptl D\subset (\ptl E\setminus L)\cap\sph_\mu$. For every $p\in\ptl D$, consider the integral curve
$\ga_{p}:I_{p}\to\hh^n$ of $J(\nuh)$, where $I_{p}$ is the maximal
interval for which $\ga_{p}$ is defined.  The trace $\ga_p(I_{p})$ is
contained in $\ptl E\setminus L$.  Such a curve is also an integral curve of
$J(X)$ in $(\hh^n)^+$ and an integral curve of $J(Y)$ in $(\hh^n)^-$,
and so it is contained in the sphere $\sph_{\mu}$.  In fact, it is
part of a characteristic curve of $\sph_{\mu}$.

Since $\ptl E\setminus L$ is foliated by integral curves of $J(\nu_H)$, it is easy to check that $\ptl E\supset \bigcup_{p\in\ptl
D}\ga_{p}(I_{p})=\sph_{\mu}\setminus (\sph_{\mu})_{0}=\sph_\mu\setminus L$. Here  $(\sph_\mu)_0$ is the singular set of the $C^2$ hypersurface $\sph_\mu$, as defined in \S~2.3.  This implies that $\sph_{\mu}\subset\ptl E$. We claim that $\bb_\mu\subset\overline{E}$. To prove this we shall show that $\bb_\mu\setminus L\subset\overline{E}$ reasoning by contradiction. If $\bb_\mu\setminus L$ is not contained in $\overline{E}$, as $\sph_\mu=\ptl\bb_\mu\subset\ptl E$, there is a point $p$ in the interior of $\bb_\mu\setminus L$ so that $p\not\in\overline{E}$. The Euclidean orthogonal projection $p'$ of $p$ over $t=0$ lies in $D\subset E$. Hence there is a point $q$ in the Euclidean segment $[p,p']\subset\intt(\bb_\mu)$ that belongs to $\ptl E\setminus L$. As $\ptl E\setminus L$ is $\hh$-regular, the perimeter of $\ptl E$ in a small ball contained in $\intt(\bb_\mu\setminus L)$ and centered at $q$ is positive, and so $|\ptl E|>|\ptl\bb_\mu|$, which contradicts our assumption that equality holds in \eqref{eq:main2}. This implies $\bb_\mu\subset\overline{E}$. As $|\bb_{\mu}|=|E|$ we obtain $\bb_{\mu}=\overline{E}$ by the normalization of $E$.

%Assume now that equality holds in \eqref{eq:main2} and that $\ptl E$
%is an $\hh$-regular hypersurface.  Then by Lemma~\ref{lem:charcurve},
%the integral curves of $J(\nu_{H})$ are contained in $\ptl E-L$.  Since
%these curves are also integral curves of $J(X)$ in $(\hh^n)^+$ or of
%$J(Y)$ on $(\hh^n)^-$, we reason as the previous case to conclude that
%$\ptl E=\sph_{\mu}$.
\end{proof}

% \begin{remark}
% Roberto Monti has recently obtained a proof of the characterization of
% equality case by assuming the isoperimetric property and the
% continuity of the horizontal unit normal.
% \end{remark}

%\noindent \emph{Added November 19th, 2009.} In a recent paper, R. Monti and D. Vittone \cite[Thm.~1,2]{mv} have proven that, given a set $E$ with finite perimeter so that the horizontal unit normal locally coincides with the restriction of a horizontal continuous vector field, the boundary $\ptl E$ is locally an $\hh$-regular surface (possibly modifying $E$ in a set of measure $0$). Using this result and Theorem~\ref{th:main} we obtain the stronger result

%\begin{theorem}
%\label{th:main2}
%Let $E\subset\hh^n$ be a normalized finite perimeter set such that $D\subset E\subset C$, where $D=D_r$, $C=C_r$, for some $r>0$. Then
%\begin{equation}
%\label{eq:main}
%|\ptl E|\ge |\ptl\bb_{\mu}|,
%\end{equation}
%where $\bb_{\mu}$ is the ball with $|\bb_{\mu}|=|E|$.  Equality holds in \eqref{eq:main} if and only if $E=\bb_{\mu}$.
%\end{theorem}

\bibliography{calibration-final}

\end{document}